\documentclass[a4paper,11pt]{article}
\usepackage{amssymb}
\usepackage{graphicx}
\usepackage{amsmath}

\input{tcilatex}

\begin{document}

\title{The number of solutions of an equation related\\
to a product of multilinear polynomials\bigskip }
\author{T. Narayaninsamy, D.-J. Mercier, J.-P. Cherdieu\medskip \medskip \\
\textit{Laboratoire A.O.C. de l'Universit\'{e} des Antilles et de la Guyane,}%
\\
97159 \textit{Pointe-\`{a}-Pitre cedex, France}.}
\date{march 1$^{\text{st}}$, 2006}
\maketitle

\parindent0cm%

\begin{center}
\textbf{Abstract}
\end{center}

\begin{quote}
We look at the number of solutions of an equation of the form $%
f_{1}f_{2}...f_{k}=a$ in a finite field, where each $f_{i}$ is a multilinear
polynomial. We use two methods to construct a solution of this problem for
the cases $a=0$, $a\neq 0$, and we generally get a semi-explicit formula. We
show that this formula can generate a more efficient algorithm\ than the
traditional algorithm which consists to make a systematic computation. We
also give explicit formulas in some special cases, and an application of our
main result to the search of the weight hierarchy of the multilinear code
with separated variables.\smallskip 
\end{quote}

\begin{center}
\textbf{Keywords}
\end{center}

\begin{quote}
Linear systems, multilinear polynomials, finite fields, separated variables,
exponential sums.\smallskip
\end{quote}

\begin{center}
\textbf{MSC Codes}
\end{center}

\begin{quote}
Primary classification: 11T06; Secondary classification: 11T23,
11T71.\bigskip \bigskip
\end{quote}

I \ INTRODUCTION\medskip

Let $\mathbb{F}_{q}$ be a finite field with $q$ elements, $\left(
m,n,k\right) \in \left( \mathbb{N}-\left\{ 0\right\} \right) ^{3}$ with $%
m\leq n,$ $\{J_{1},J_{2},...,J_{m}\}$ a partition of $\left\{
1,...,n\right\} $, $\left( f,a\right) \in \mathbb{F}_{q}\left[ X_{1},..,X_{n}%
\right] \times \mathbb{F}_{q}$, $A$ an $k\times m$ matrix with entries in $%
\mathbb{F}_{q}$. We assume that $f$ satisfies the following:\medskip

$\left( 1.1\right) $\ \ \ \ \ \ \ \ \ \ \ \ \ \ \ \ \ \ \ \ \ \ \ \ \ \ \ \
\ \ \ \ \ \ \ \ \ $f=f_{1}f_{2}...f_{k},$

\quad where :

$\left( 1.2\right) $\ \ \ \ \ \ \ \ \ \ \ \ \ \ \ \ \ \ \ \ \ \ \ \ \ \ \ \
\ \ \ \ \ $f_{i}\left( X_{1},..,X_{n}\right)
=\sum\limits_{j=1}^{m}a_{ij}\prod\limits_{\tau \in J_{j}}X_{\tau }$ , \ $%
A=\left( a_{ij}\right) _{\substack{ 1\leq i\leq k  \\ 1\leq j\leq m}}$%
.\medskip

The polynomials $f_{i}$ are called multilinear polynomials, ($\left[ 3\right]
$ p.~32)$.$ In this present paper, the problem of determining the number of
solutions in $\mathbb{F}_{q}^{n}$ of the following equation is
considered:\medskip

$\left( 1.3\right) $ \ \ \ \ \ \ \ \ \ \ \ \ \ \ \ \ \ \ \ \ \ \ \ \ \ \ \ \
\ \ \ \ \ $f\left( X_{1},..,X_{n}\right) =a.$\medskip

The solution of this problem is known only for special cases of $\left(
1.1\right) $ ($\left[ 1\right] $, $\left[ 2\right] $, $\left[ 3\right] $).
We determine for the general case a semi-explicit formula for the solution
of this problem. We obtain an explicit formula in special cases. When $a\neq
0,$ the considered method is based on a generalization of the following
formula:\medskip

$\left( 1.4\right) $ \ \ \ \ \ \ \ \ \ \ \ \ \ \ \ \ \ $N\left( h\left(
X\right) g\left( X\right) =a\right) =\sum\limits_{u\in \mathbb{F}_{q}^{\ast
}}\left( N\left( h\left( X\right) =u\text{ \& }g\left( X\right) =\frac{a}{u}%
\right) \right) ,$\ 

where $N\left( h\left( X\right) =u\text{ \& }g\left( X\right) =a/u\right) $
denotes the number of common solutions of$\ $the equations $h\left( X\right)
=u$ and $g\left( X\right) =a/u.$\medskip

When $a=0$,\ the considered method is based on the so-called '' inclusion -
exclusion principle''. Note that if $m\geq 2$ and if $a_{ij}$ $\neq 0$ for
at least two entries $j,$ the polynomials of the form $\left( 1.2\right) $
are strictly contained in the class of irreducible polynomials. Then, it
could be interesting in a futur study to determine a sufficient condition\
to have the form $\left( 1.1\right) $. When this form exists, there exist
algorithms to find these irreducible factors $\left[ 4\right] $. In our
study, we suppose that the associated irreducible factors\ to $f$ are known
and our aim is to determine an expression for the number of solutions of the
equation~$\left( 1.3\right) .$\medskip

The organisation of this paper is the following: In section II, we continue
the presentation of the considered notations, we recall basic results
concerning linear systems and exponential sums. The ''inclusion - exclusion
principle'' is also recalled. We compute formulas for the numbers of
solutions of special systems of equations in Section III. These formulas
will be useful in the treatment of the considered problem. The special case $%
k\in \left\{ 1,2\right\} $\ is considered in this Section.\medskip

The main result for the general case is Proposition IV.1. In section V, we
apply the obtained result to two ''numerical'' examples such that $k\in
\left\{ 2,3\right\} $, $a_{ij}\in \left\{ 0,1\right\} ,$ $m\in \left\{
3,4\right\} $, $n\in \left\{ 7,8\right\} .$ We show from these examples that
the obtained formula can generate a clearly more efficient algorithm than
the ''traditional'' one which consists to make a systematic computation of
all the solutions of this equation.\medskip

In section VI we show that the research of the number of zeros of products $%
f=f_{1}f_{2}...f_{k}$ of the form $\left( 1.1\right) $ is important when we
compute the weight hierarchy of a multilinear code. To conclude in Section
VII, we give three examples of problems that can be investigated from the
present study.

\ 

\ 

\ 

II \ NOTATIONS, BASIC RESULTS\medskip

Recall that $m\leq n$. We assume:

\ 

$\left( 2.1\right) $\ \ \ \ \ \ \ \ \ $k\leq m,$

\ 

$\left( 2.2\right) $\ \ \ \ \ \ \ \ $A$ is an $k\times m$ matrix of rank of $%
k.$

\ 

Let \ \ \ \ \ $1\leq i_{1}<i_{2}<...<i_{l}\leq k,$ \ $\left(
b_{1},...,b_{l}\right) \in \mathbb{F}_{q}^{l}.$

\ 

The system:

\ 

$\left( 2.3\right) $ $\ \ \ \ \ \ \ \ \ \ \ $\ $\ \left\{ 
\begin{array}{c}
\ f_{i_{1}}\left( X_{1},..,X_{n}\right) =b_{1}, \\ 
\\ 
\ f_{i_{2}}\left( X_{1},..,X_{n}\right) =b_{2}, \\ 
. \\ 
. \\ 
f_{i_{l}}\left( X_{1},..,X_{n}\right) =b_{l},%
\end{array}
\right. $

\ 

is equivalent to:

\ 

\ \ \ \ \ \ \ \ \ \ \ $A_{\left[ i_{1},i_{2},...,i_{l}\right] }\left(
\prod\limits_{\tau \in J_{1}}X_{\tau },\prod\limits_{\tau \in J_{2}}X_{\tau
},....,\prod\limits_{\tau \in J_{m}}X_{\tau }\right) ^{T}=\left(
b_{1},...,b_{l}\right) ^{T},$

\ 

where

\ 

\ \ \ \ \ \ \ \ \ \ $A_{\left[ i_{1},i_{2},...,i_{l}\right] }$ =$\left(
a_{i_{t,}j}\right) _{\substack{ 1\leq t\leq l  \\ 1\leq j\leq m}}.$

\ 

We can write:

\ 

$\left( 2.4a\right) $ $\ \ \ \ \ \ \ \ \ \ B_{\left[ i_{1},i_{2},...,i_{l}%
\right] }\left( \prod\limits_{\tau \in J_{\nu _{1}}}X_{\tau
},\prod\limits_{\tau \in J_{\nu _{2}}}X_{\tau },....,\prod\limits_{\tau \in
J_{\nu _{l}}}X_{\tau }\right) ^{T}=\left( b_{1},...,b_{l}\right) ^{T}-$

\ \ 

$\ \ \ \ \ \ \ \ \ \ \ \ \ \ \ \ \ \ \ \ \ \ \ \ \ \ \ \ \ \ \ \ \ \ \ \ \ \
\ \ \ \ \ \ \ \ \ \ C_{\left[ i_{1},i_{2},...,i_{l}\right] }\left(
\prod\limits_{\tau \in J_{\nu _{1}^{,}}}X_{\tau },\prod\limits_{\tau \in
J_{\nu _{2}^{,}}}X_{\tau },....,\prod\limits_{\tau \in J_{\nu
_{m-l}^{,}}}X_{\tau }\right) ^{T},$

\ 

where

\ 

$\left( 2.4b\right) $\ \ \ \ \ \ \ \ \ \ \ \ \ $B_{\left[
i_{1},i_{2},...,i_{l}\right] }=\left( a_{i_{t,}j}\right) _{\substack{ 1\leq
t\leq l  \\ j\in \left\{ \nu _{1,...,}\nu _{l}\right\} \subset \left\{
1,...,m\right\} }}$

\ 

\ \ \ \ \ \ \ \ \ \ \ \ \ \ \ \ \ \ \ \ is an invertible $l\times l$
submatrix of $A_{\left[ i_{1},i_{2},...,i_{l}\right] }$,

\ \ 

$\left( 2.4c\right) $\ \ \ \ \ \ \ \ \ $C_{\left[ i_{1},i_{2},...,i_{l}%
\right] }=\left( a_{i_{t,}j}\right) _{\substack{ 1\leq t\leq l  \\ j\in
\left\{ 1,...,m\right\} -\left\{ \nu _{1,...,}\nu _{l}\right\} }}$

\ \ \ 

\ \ \ \ \ \ \ \ \ \ \ \ \ \ \ \ \ \ \ is an $l\times \left( m-l\right) $ \
submatrix of $A,$

\ 

$\left( 2.4d\right) $ \ \ \ \ \ \ \ \ $\left\{ \nu _{1}^{,},...,\nu
_{m-l}^{,}\right\} =\left\{ 1,...,m\right\} -\left\{ \nu _{1,...,}\nu
_{l}\right\} .$

\ 

We will use the following notations:

\ \ \ \ \ \ \ \ \ \ \ \ \ \ \ \ \ \ \ \ \ \ 

$\left( 2.5a\right) $ \ \ \ \ \ \ \ \ $\left\{ 
\begin{array}{c}
B_{\left[ i_{1},i_{2},...,i_{l}\right] }^{-1}=\left( a_{i_{t,}j}^{^{\prime
}}\right) _{\substack{ t\in \left\{ 1,...,l\right\}  \\ j\in \left\{ \nu
_{1,...,}\nu _{l}\right\} \subset \left\{ 1,...,m\right\} }}, \\ 
\\ 
B_{\left[ i_{1},i_{2},...,i_{l}\right] }^{-1}C_{\left[ i_{1},i_{2},...,i_{l}%
\right] }=\left( \sigma _{\left[ i_{1},i_{2},...,i_{l}\right] }^{t,j}\right) 
_{\substack{ 1\leq t\leq l  \\ 1\leq j\leq m-l}}, \\ 
\\ 
I_{\nu ^{\prime }}=J_{\nu _{1}^{,}}\cup J_{\nu _{2}^{,}}\cup ..\cup J_{\nu
_{m-l}^{,}}.%
\end{array}
\right. $

\ 

We denote by

\ 

$\left( 2.5b\right) \ \ \ \ \ \ \ \ \ \ \ \ \ \ \ \ \ \ \ \ \ \ N\left(
f_{i_{1}},....,f_{i_{l}},b_{1},...,b_{l}\right) ,$

\ 

the number of solutions in $\mathbb{F}_{q}^{n}$ of the system of equations $%
\left( 2.3\right) $. The number of solutions in $\mathbb{F}_{q}^{n}$ of the
equation $\left( 1.3\right) $ is denoted by

\ 

$\left( 2.5c\right) \ \ \ \ \ \ \ \ \ \ \ \ \ \ \ \ \ \ \ \ \ \ N\left(
f,a\right) .$\ \ 

\ 

To treat the considered problem, we use a method based on the following
relations:

\ 

$\left( 2.6a\right) $\ \ \ \ \ \ \ \ \ $N\left( f,a\right) =\sum\limits 
_{\substack{ \left( a_{1},...,a_{k}\right) \in \mathbb{F}_{q}^{k}  \\ %
a_{1}....a_{k}=a}}N\left( f_{1},...,f_{k},a_{1},...,a_{k}\right) $ , \ 

\ 

\ \ \ \ \ \ \ \ \ \ \ \ \ \ \ \ \ \ \ \ when $a\neq $ $0,$

\ 

$\left( 2.6b\right) \ \ \ \ \ \ \ \ \ N\left( f,0\right) =N\left(
f_{1},0\right) +N\left( f_{2},0\right) +...+N\left( f_{k},0\right)
-\sum\limits_{1\leq i_{1}<i_{2}\leq k}N\left( f_{i_{1}},f_{i_{2}},0,0\right) 
$

$\bigskip $

$\ \ \ \ \ +\sum\limits_{1\leq i_{1}<i_{2}<i_{3}\leq k}N\left(
f_{i_{1}},f_{i_{2}},f_{i_{3}},0,0,0\right) +....+\left( -1\right)
^{j+1}\sum\limits_{1\leq i_{1}<i_{2}..<i_{j}\leq k}N\left(
f_{i_{1}},f_{i_{2}},...,f_{i_{j}},0,0,...,0\right) $

\ 

\bigskip\ $\ \ \ \ \ \ \ \ \ \ \ \ \ \ +...+\left( -1\right)
^{k+1}N(f_{1},f_{2},...,f_{k},0,0,...,0)$\ $).$

\ 

The considered method in $\left[ 2\right] $ to solve this problem in the
particular case $k=1,$ $a=0$, is based on exponential sums. We shall use a
similar method in Section III to solve the more general case $k=1$.

For all $u\in \mathbb{F}_{q}$, the following function is an additive
character on $\mathbb{F}_{q}$

\ 

$\left( 2.7\right) $ \ \ \ \ \ \ \ \ \ \ \ \ \ \ \ \ $\Psi _{u}$\ $\left(
v\right) =\exp \left( \frac{2i\pi }{p}Tr_{\mathbb{F}_{q}}\left( uv\right)
\right) ,$

\ 

where $Tr_{\mathbb{F}_{q}}\left( uv\right) $\ is the absolute trace of $uv.$

\ 

We will use the following basic results:

\ 

$\left( 2.8\right) $ \ \ \ \ \ \ \ \ \ \ \ $\sum\limits_{u\in \mathbb{F}%
_{q}}\Psi _{u}$\ $\left( v\right) =\left\{ 
\begin{array}{c}
0\text{ \ if \ }v\neq 0, \\ 
q\text{ \ if \ }v=0.%
\end{array}
\right. $\ 

$\left( 2.9\right) \ \ \ \ \ \ \ \ \ \ \ \sum\limits_{v\in \mathbb{F}%
_{q}}\Psi _{u}$\ $\left( v\right) =\left\{ 
\begin{array}{c}
0\text{ \ if \ }u\neq 0, \\ 
q\text{ \ if \ }u=0.%
\end{array}
\right. $\ \ 

\ 

We also need the following result given in $\left[ 2\right] $:

\ 

\textbf{PROPOSITION II-1}

\ 

$Let$ \ $(\alpha ,u,d)\in (\mathbb{F}_{q}-\left\{ 0\right\} )^{2}\times (%
\mathbb{N}-\left\{ 0\right\} ),$ $then$:

\ 

$\left( 2.10\right) $ \ \ \ \ \ \ \ \ \ \ $\sum\limits_{\left(
X_{1},X_{2},..,X_{d}\right) \in \mathbb{F}_{q}^{d}}$\ \ $\Psi _{u}$\ $\left(
\alpha \prod\limits_{i\in \left\{ 1,..,d\right\} }X_{i}\right)
=q^{d}-q\left( q-1\right) ^{d-1}.$

\ 

We conclude this present Section by recalling the wellknown map $\left[ 1%
\right] $:

\ 

$\left( 2.11\right) $\ \ \ \ \ \ \ \ \ \ \ \ \ \ \ \ \ \ \ \ \ \ \ \ \ \ \ \
\ \ $\kappa $ : $\mathbb{F}_{q}\rightarrow \left\{ -1,q-1\right\} ,$

\ \ \ \ \ \ \ \ \ \ \ \ \ \ \ \ \ \ \ \ \ \ \ \ \ \ \ \ \ \ $\ \ \ \ \ \ \
\kappa \left( X\right) =\left\{ 
\begin{array}{l}
-1\text{ \ \ \ \ if \ }X\neq 0, \\ 
q-1\text{ \ if \ }X=0.%
\end{array}
\right. \bigskip \bigskip $

III \ NUMBER OF\ SOLUTIONS OF SPECIAL\ SYSTEMS\medskip

In the following Proposition, we determine the solution of the problem for
the special case $k=1$.\ 

\ 

\textbf{PROPOSITION III-1}

\thinspace

$Assume$ $that$ $k=1$ $in$ $\left( 1.1\right) ,$ $then$:

\ 

$\left( 3.1\right) $\ \ \ \ \ \ \ \ \ \ $N\left( f,a\right) =q^{n-1}+\kappa
\left( a\right) q^{n-1}\left[ \prod\limits_{\left\{ l/a_{1l}\neq 0\right\}
}\left( 1-\left( \frac{q-1}{q}\right) ^{\left| J_{l}\right| -1}\right) %
\right] ,$

\ 

\ \ \ \ \ \ \ \ \ \ \ $\left( Recall\text{ \ }that\text{ }N\left( f,a\right) 
\text{ }is\text{ }the\text{ }number\text{ }of\text{ \ }solutions\text{ }of%
\text{ \ }\left( 1.3\right) \right) .$

\thinspace

\textbf{Proof}

\thinspace

First, $\left( 2.8\right) $ and $\left( 2.9\right) $ imply that:

\ 

\ \ \ $N\left( f,a\right) =\frac{1}{q}\left( \sum\limits_{\left(
X_{1},X_{2},..,X_{n}\right) \in \mathbb{F}_{q}^{n}}\sum\limits_{u\in \mathbb{%
F}_{q}}\Psi _{u}\ (f\left( X_{1},X_{2},..,X_{n}\right) -a)\right) ,$

\ 

thus:

\ \ \ 

\ \ \ \ $N\left( f,a\right) =$\ \ \ \ \ $\frac{1}{q}\left( \sum\limits_{u\in 
\mathbb{F}_{q}}\sum\limits_{\left( X_{1},X_{2},..,X_{n}\right) \in \mathbb{F}%
_{q}^{n}}\Psi _{u}\ (f\left( X_{1},X_{2},..,X_{n}\right) -a)\right) $

\ 

\ \ \ \ \ \ \ \ \ \ \ \ \ \ \ $=\frac{1}{q}\left( \sum\limits_{u\in \mathbb{F%
}_{q}}(\Psi _{u}\ (-a)\sum\limits_{\left( X_{1},X_{2},..,X_{n}\right) \in 
\mathbb{F}_{q}^{n}}\Psi _{u}\ (f\left( X_{1},X_{2},..,X_{n}\right) ))\right) 
$

\ \ \ \ \ \ \ \ \ \ \ \ \ \ \ $=$\ $\frac{1}{q}\left( \sum\limits_{u\in 
\mathbb{F}_{q}}(\Psi _{u}\ (-a)\sum\limits_{\left(
X_{1},X_{2},..,X_{n}\right) \in \mathbb{F}_{q}^{n}}\Psi
_{u}(\sum\limits_{j=1}^{m}a_{1j}\prod\limits_{\tau \in J_{j}}X_{\tau
}))\right) $

\ 

\ \ \ \ \ \ \ \ \ \ \ \ \ \ \ $=\frac{1}{q}$($\sum\limits_{u\in \mathbb{F}%
_{q}}[\Psi _{u}\ (-a)\prod\limits_{\left\{ k/a_{1k}\neq 0\right\}
}\sum\limits_{\left( X_{1},X_{2},..,X_{\left| J_{k}\right| }\right) \in 
\mathbb{F}_{q}^{\left| J_{k}\right| }}\Psi _{u}(a_{1k}\prod\limits_{\tau \in
J_{k}}X_{\tau })]).$

\ 

After, from Proposition II-1 we get:

\ 

\ \ \ \ $N\left( f_{1},a\right) =q^{n-1}+q^{n-1-\sum\limits_{\left\{
k/a_{1k}\neq 0\right\} }\left| J_{k}\right| }\left( \sum\limits_{u\in 
\mathbb{F}_{q}^{\ast }}[\Psi _{u}\ (-a)\prod\limits_{\left\{ k/a_{1k}\neq
0\right\} }\left( q^{\left| J_{k}\right| }-q\left( q-1\right) ^{\left|
J_{k}\right| -1}\right) \right) $

\ 

\ \ \ \ \ \ \ \ \ \ \ \ \ \ \ \ $=q^{n-1}+q^{n-1}[\prod\limits_{\left\{
k/a_{1k}\neq 0\right\} }\left( 1-\left( \frac{q-1}{q}\right) ^{\left|
J_{k}\right| -1}\right) ][\sum\limits_{u\in \mathbb{F}_{q}^{\ast }}\Psi
_{u}\ (-a)].$

\ 

From $\left( 2.8\right) \ $and $\left( 2.11\right) $ we deduce $\sum_{u\in 
\mathbb{F}_{q}^{\ast }}\Psi _{u}(-a)=\kappa \left( a\right) $, thus

\ 

\ \ \ \ \ \ \ \ \ \ \ \ \ \ \ \ \ \ $N\left( f,\text{ }a\right)
=q^{n-1}+\kappa \left( a\right) q^{n-1}\prod\limits_{\left\{ l/a_{1l}\neq
0\right\} }\left( 1-\left( \frac{q-1}{q}\right) ^{\left| J_{l}\right|
-1}\right) $,

\ \ \ \ 

which is the desired\ result.$\blacksquare $

\ 

The following Corollary which is a consequence of Proposition III-1, will be
useful for instance in the proof of Proposition III-3.\bigskip \bigskip\ 

\textbf{COROLLARY III-2}

\thinspace

$Let$ $\left( w,\alpha ,d\right) \in \mathbb{F}_{q}\times \left( \mathbb{F}%
_{q}-\left\{ 0\right\} \right) \times \left( \mathbb{N}-\left\{ 0\right\}
\right) ,$ $then$ $the$ $number$ $N\left( \alpha \prod\limits_{\tau
=1}^{d}X_{\tau },w\right) $

$of$ $solutions$ $in$ $\mathbb{F}_{q}^{d}$ $of$ $the$\ $equation$:

\ 

\ \ \ \ \ \ \ \ \ \ \ \ \ \ \ \ \ \ \ \ \ \ \ \ \ \ \ \ \ \ \ $\alpha
\prod\limits_{\tau =1}^{d}X_{\tau }=w$

$is$

\ \ \ \ \ \ \ \ \ \ \ \ \ \ \ \ \ $N\left( \alpha \prod\limits_{\tau
=1}^{d}X_{\tau },w\right) =$\ $q^{d-1}+\kappa \left( w\right) q^{d-1}\left(
1-\left( \frac{q-1}{q}\right) ^{d-1}\right) .$

\bigskip

\textbf{PROPOSITION III-3}

$\,$

$Let$ \ \ \ $1\leq i_{1}<i_{2}<...<i_{l}\leq k,$ $\left(
b_{1},...,b_{l}\right) \in \mathbb{F}_{q}^{l},$ $then$ \ 

\ 

$\left( 3.2\right) $ $\ \ \ \ \ \ N\left(
f_{i_{1}},f_{i_{2}},...,f_{i_{l}},b_{1},b_{2},...,b_{l}\right) =$

\ 

\ 

$\sum\limits_{\left( X_{i}\right) _{i\in I_{\nu ^{\prime }}}\in \mathbb{F}%
_{q}^{\left| I_{\nu ^{\prime }}\right| }}\prod\limits_{t=1}^{l}\left(
q^{\left| J_{\nu _{t}}\right| -1}+\kappa \left(
\sum\limits_{j=1}^{l}a_{i_{t,}j}^{^{\prime
}}b_{j}-\sum\limits_{j=1}^{m-l}\sigma _{\left[ i_{1},i_{2},...,i_{l}\right]
}^{t,j}\prod\limits_{\tau \in J_{\nu _{j}^{,}}}X_{\tau }\right) q^{\left|
J_{\nu _{t}}\right| -1}\left( 1-\left( \frac{q-1}{q}\right) ^{\left| J_{\nu
_{t}}\right| -1}\right) \right) ,$

\ 

\ \ \ \ \ \ \ \ \ \ \ \ \ \ \ \ \ \ $\left( I_{\nu ^{\prime }},\text{ }%
J_{\nu _{t}^{,}}\text{ , }a_{i_{t,}j}^{^{\prime }}\text{ , }\sigma _{\left[
i_{1},i_{2},...,i_{l}\right] }^{t,j}\text{ }are\text{ }defined\text{ }in%
\text{ }\left( 2.4b\right) ,\left( 2.4c\right) ,\left( 2.4d\right) ,\left(
2.5a\right) \right) .$

\ \ \ \ \ \ \ \ \ \ \ 

\textbf{Proof}

\thinspace

First recall that the matrix $A$ is of rank $k,$ then\ the existence of the $%
l\times l$ invertible matrix \ $B_{\left[ i_{1},i_{2},...,i_{l}\right] }$
introduced in $\left( 2.4a\right) ,$ $\left( 2.4b\right) $ is justified.
Consequently, the \ equality $\left( 2.4\right) $ is equivalent to:

\ 

\ \ $\left( \prod\limits_{\tau \in J_{\nu _{1}}}X_{\tau },\prod\limits_{\tau
\in J_{\nu _{2}}}X_{\tau },....,\prod\limits_{\tau \in J_{\nu _{l}}}X_{\tau
}\right) ^{T}=B_{\left[ i_{1},i_{2},...,i_{l}\right] }^{-1}\left(
b_{1},...,b_{l}\right) ^{T}-$

\ \ 

$\ \ \ \ \ \ \ \ \ \ \ \ \ \ \ \ \ \ \ \ \ \ \ \ \ \ \ \ \ \ \ \ \ \ \ \ B_{%
\left[ i_{1},i_{2},...,i_{l}\right] }^{-1}C_{\left[ i_{1},i_{2},...,i_{l}%
\right] }\left( \prod\limits_{\tau \in J_{\nu _{1}^{,}}}X_{\tau
},\prod\limits_{\tau \in J_{\nu _{2}^{,}}}X_{\tau },....,\prod\limits_{\tau
\in J_{\nu _{m-l}^{,}}}X_{\tau }\right) ^{T},$

\ 

\ 

\ \ \ with \ \ \ \ \ $\left\{ \nu _{1,...,}\nu _{l}\right\} \cup \left\{ \nu
_{1}^{,},...,\nu _{m-l}^{,}\right\} =\left\{ 1,...,m\right\} ,$ $\left( 
\text{see }\left( 2.4d\right) \right) .$

\ 

It follows that:

$\left( 3.2\right) $\ \ \ \ \ \ \ \ \ \ \ \ \ \ \ \ $\prod\limits_{\tau \in
J_{\nu _{t}}}X_{\tau }=\sum\limits_{j=1}^{l}a_{i_{t,}j}^{^{\prime
}}b_{j}-\sum\limits_{j=1}^{m-l}\sigma _{\left[ i_{1},i_{2},...,i_{l}\right]
}^{t,j}\prod\limits_{\tau \in J_{\nu _{j}^{,}}}X_{\tau }$ \ , \ \ $t=1,...,l$%
,

\ 

applying Corollary III-2\ to $\left( 3.2\right) ,$ we get:

\ 

$N\left( f_{i_{1}},f_{i_{2}},...,f_{i_{j}},b_{1},b_{2},...,b_{l}\right) =$

\ 

$\sum\limits_{\left( X_{i}\right) _{i\in I_{\nu ^{\prime }}}\in \mathbb{F}%
_{q}^{\left| I_{\nu ^{\prime }}\right| }}\prod\limits_{t=1}^{l}\left(
q^{\left| J_{\nu _{t}}\right| -1}+\kappa \left(
\sum\limits_{j=1}^{l}a_{i_{tt,}j}^{^{\prime
}}b_{j}-\sum\limits_{j=1}^{m-l}\sigma _{\left[ i_{1},i_{2},...,i_{l}\right]
}^{t,j}\prod\limits_{\tau \in J_{\nu _{j}^{,}}}X_{\tau }\right) q^{\left|
J_{\nu _{t}}\right| -1}\left( 1-\left( \frac{q-1}{q}\right) ^{\left| J_{\nu
_{t}}\right| -1}\right) \right) ,$

\ 

which is the announced result.$\blacksquare $

\ 

Now in the following Proposition, we give a semi-explicit formula for the
number of solutions\ of $\left( 1.3\right) $ for the case $k=2.$

\ 

\textbf{PROPOSITION III-4}

$\;$

$Assume$ $that$ $k=2\ in$ $\left( 1.1\right) $. $With$ $the$ $notations$ $%
\left( 2.5a\right) $, $\left( 2.5b\right) \ and$ $\left( 2.5c\right) ,\;we$ $%
have$:

\ 

$\ \ \ \ \ N\left( f,0\right) =2q^{n-1}+$

$\ \ \ \ \ \ \ \ \ \ \ \ \ \ \ \ \ \ \ \ \ \ \ \ \ \ \ \ \left( q-1\right)
q^{n-1}\left( \prod\limits_{\left\{ l/a_{1l}\neq 0\right\} }\left( 1-\left( 
\frac{q-1}{q}\right) ^{\left| J_{l}\right| -1}\right) +\prod\limits_{\left\{
l/a_{2l}\neq 0\right\} }\left( 1-\left( \frac{q-1}{q}\right) ^{\left|
J_{l}\right| -1}\right) \right) $

\ \ \ 

\ \ \ \ \ \ \ \ \ \ \ \ $-\sum\limits_{\left( X_{i}\right) _{i\in I_{\nu
^{\prime }}}\in \mathbb{F}_{q}^{\left| I_{\nu ^{\prime }}\right|
}}\prod\limits_{t=1}^{2}\left( q^{\left| J_{\nu _{t}}\right| -1}+\kappa
\left( \sum\limits_{j=1}^{m-2}\sigma _{\left[ 1,2\right] }^{t,j}\prod%
\limits_{\tau \in J_{\nu _{j}^{,}}}X_{\tau }\right) q^{\left| J_{\nu
_{t}}\right| -1}\left( 1-\left( \frac{q-1}{q}\right) ^{\left| J_{\nu
_{t}}\right| -1}\right) \right) ,$

\ \ \ 

\ \ 

\ \ \ \ \ $N\left( f,a\right) =$

\ 

$\sum\limits_{u\in \mathbb{F}_{q}^{\ast }}\sum\limits_{\left( X_{i}\right)
_{i\in I_{\nu ^{\prime }}}\in \mathbb{F}_{q}^{\left| I_{\nu ^{\prime
}}\right| }}\prod\limits_{t=1}^{2}q^{\left| J_{\nu t}\right| -1}\left(
1+\kappa \left( a_{1,1}^{^{\prime }}u+\frac{a_{1,2}^{^{\prime }}}{u}%
a-\sum\limits_{j=1}^{m-2}\sigma _{\left[ 1,2\right] }^{t,j}\prod\limits_{%
\tau \in J_{\nu _{j}^{,}}}X_{\tau }\right) \left( 1-\left( \frac{q-1}{q}%
\right) ^{\left| J_{\nu _{t}}\right| -1}\right) \right) ,$

\ 

\ \ \ \ \ \ \ $if$ \ \ \ $a\neq 0.$

\ \ 

\textbf{Proof}

\thinspace

First consider the case $a=0.$ Then, from $\left( 2.6b\right) $ we have:

\ 

$\ \ \ \ \ \ \ \ \ \ \ N\left( f,0\right) =N\left( f_{1},0\right) +N\left(
f_{2},0\right) -N\left( f_{1},f_{2},0,0\right) .$

\ 

Using $\left( 3.1\right) \ $ we get:

\ 

$\left( 3.3\right) $\ \ $\ \ \ \ \ N\left( f_{1},0\right) +N\left(
f_{2},0\right) =$

\ 

$\ \ q^{n-1}+\kappa \left( 0\right) q^{n-1}\left[ \prod\limits_{\left\{
l/a_{1l}\neq 0\right\} }\left( 1-\left( \frac{q-1}{q}\right) ^{\left|
J_{l}\right| -1}\right) \right] +q^{n-1}+\kappa \left( 0\right) q^{n-1}\left[
\prod\limits_{\left\{ l/a_{2l}\neq 0\right\} }\left( 1-\left( \frac{q-1}{q}%
\right) ^{\left| J_{l}\right| -1}\right) \right] $

\ 

$\ \ =2q^{n-1}+\left( q-1\right) q^{n-1}\left[ \prod\limits_{\left\{
l/a_{1l}\neq 0\right\} }\left( 1-\left( \frac{q-1}{q}\right) ^{\left|
J_{l}\right| -1}\right) +\prod\limits_{\left\{ l/a_{2l}\neq 0\right\}
}\left( 1-\left( \frac{q-1}{q}\right) ^{\left| J_{l}\right| -1}\right) %
\right] .$

\ 

After, using Proposition III-3 in the case $\left( l,b_{1},b_{2}\right)
=\left( 2,0,0\right) $, we can write:

\ 

$N\left( f_{1},f_{2},0,0\right) =\sum\limits_{\left( X_{i}\right) _{i\in
I_{\nu ^{\prime }}}\in F_{q}^{\left| I_{\nu ^{\prime }}\right|
}}\prod\limits_{t=1}^{2}\left( q^{\left| J_{\nu _{t}}\right| -1}+\kappa
\left( \sum\limits_{j=1}^{m-2}\sigma _{\left[ 1,2\right] }^{t,j}\prod%
\limits_{\tau \in J_{\nu _{j}^{,}}}X_{\tau }\right) q^{\left| J_{\nu
_{t}}\right| -1}\left( 1-\left( \frac{q-1}{q}\right) ^{\left| J_{\nu
_{t}}\right| -1}\right) \right) ,$

\ 

the desired result follows in the case $a=0.$

\ 

Now assume that $a\neq 0.$ From $\left( 2.6a\right) $ we have: 
\begin{equation*}
N\left( f,a\right) =\sum\limits_{u\in \mathbb{F}_{q}^{\ast }}N\left(
f_{1},f_{2},u,\frac{a}{u}\right) .
\end{equation*}
Let $u\in \mathbb{F}_{q}^{\ast },$ by Proposition III-3 in the case $\left(
l,b_{1},b_{2}\right) =\left( 2,u,\frac{a}{u}\right) $ we get:

\ 

$N\left( f_{1},f_{2},u,\frac{a}{u}\right) =$

\ 

$\sum\limits_{\left( X_{i}\right) _{i\in I_{\nu ^{\prime }}}\in \mathbb{F}%
_{q}^{\left| I_{\nu ^{\prime }}\right| }}\prod\limits_{t=1}^{2}q^{\left|
J_{\nu _{t}}\right| -1}\left( 1+\kappa \left( a_{1,1}^{^{\prime }}u+\frac{%
a_{1,2}^{^{\prime }}}{u}a-\sum\limits_{j=1}^{m-2}\sigma _{\left[ 1,2\right]
}^{t,j}\prod\limits_{\tau \in J_{\nu _{j}^{,}}}X_{\tau }\right) \left(
1-\left( \frac{q-1}{q}\right) ^{\left| J_{\nu _{t}}\right| -1}\right)
\right) ,$

\ 

consequently, we easily get the announced\ result in the case $a\neq 0$. $%
\blacksquare $

$\ $\ \ \ \ \ \ 

\ \ \ 

\ \ IV \ THE GENERAL PROBLEM

\ 

Here, we consider the equation $\left( 1.3\right) $, namely: $f\left(
X_{1},..,X_{n}\right) =a$, where $f=f_{1}f_{2}...f_{k}$ and $f_{i}\left(
X_{1},..,X_{n}\right) =\sum_{j=1}^{m}a_{ij}\prod_{\tau \in J_{j}}X_{\tau }$,
and we assume that the rank of the matrix $A=\left( a_{ij}\right) _{1\leq
i\leq k,1\leq j\leq m}$ is $k$. In the general case, the number of solutions 
$N\left( f,a\right) $ of $f\left( X_{1},..,X_{n}\right) =a$ is given by next
Proposition:

\ 

\textbf{PROPOSITION IV-1}

$\;$

$With$ $the$ $notations$ $\left( 2.5a\right) $,

\ \ \ 

\ \ \ \ \ \ \ $N\left( f,0\right) $ $=kq^{n-1}+\left( q-1\right)
q^{n-1}\left( \sum\limits_{i=1}^{k}\prod\limits_{\left\{ l/a_{il}\neq
0\right\} }\left( 1-\left( \frac{q-1}{q}\right) ^{\left| J_{l}\right|
-1}\right) \right) +$

\ 

$\sum\limits_{l=2}^{k}\left( -1\right) ^{l+1}\sum\limits_{1\leq
i_{1}<i_{2}<..<i_{l}\leq k}\sum\limits_{\left( X_{i}\right) _{i\in I_{\nu
^{\prime }}}\in \mathbb{F}_{q}^{\left| I_{\nu ^{\prime }}\right|
}}\prod\limits_{t=1}^{l}q^{\left| J_{\nu _{t}}\right| -1}\left( 1+\kappa
\left( \sum\limits_{j=1}^{m-l}\sigma _{\left[ i_{1},i_{2},...,i_{l}\right]
}^{t,j}\prod\limits_{\tau \in J_{\nu _{j}^{,}}}X_{\tau }\right) \left(
1-\left( \frac{q-1}{q}\right) ^{\left| J_{\nu _{t}}\right| -1}\right)
\right) ,$

\ 

$\ \ \ \ \ \ \ N\left( f,a\right) =$

\ 

$\sum\limits_{\substack{ \left( a_{1},...,a_{k}\right) \in \mathbb{F}%
_{q}^{k}  \\ a_{1}a_{2}..a_{k}=a}}\sum\limits_{\left( X_{i}\right) _{i\in
I_{\nu ^{\prime }}}\in \mathbb{F}_{q}^{\left| I_{\nu ^{\prime }}\right|
}}\prod\limits_{t=1}^{k}q^{\left| J_{\nu _{t}}\right| -1}\left( 1+\kappa
\left( \sum\limits_{j=1}^{k}a_{t,j}^{^{\prime
}}a_{j}-\sum\limits_{j=1}^{m-k}\sigma _{\left[ 1,...,k\right]
}^{t,j}\prod\limits_{\tau \in J_{\nu _{j}^{,}}}X_{\tau }\right) \left(
1-\left( \frac{q-1}{q}\right) ^{\left| J_{\nu _{t}}\right| -1}\right)
\right) ,$

$\ if$ \ $a\neq 0.$

\ 

\textbf{Proof}

\thinspace

First consider the case $a=0.$ From Proposition III-1 and for all $i\in
\left\{ 1,..,k\right\} $:

\ 

\ \ \ \ \ \ \ \ $N\left( f_{i},0\right) =q^{n-1}+\left( q-1\right) q^{n-1}%
\left[ \prod\limits_{\left\{ l/a_{il}\neq 0\right\} }\left( 1-\left( \frac{%
q-1}{q}\right) ^{\left| J_{l}\right| -1}\right) \right] ,$

thus:

\ 

\ \ \ \ \ \ $\sum\limits_{i=1}^{k}$\ $N\left( f_{i},0\right)
=kq^{n-1}+\left( q-1\right) q^{n-1}\ \sum\limits_{i=1}^{k}\left[
\prod\limits_{\left\{ l/a_{il}\neq 0\right\} }\left( 1-\left( \frac{q-1}{q}%
\right) ^{\left| J_{l}\right| -1}\right) \right] .$

\ 

Now let \ $1\leq i_{1}<i_{2}<...<i_{l}\leq k,$ Proposition III-3 gives:

\ 

\ \ \ \ \ \ \ $\ \ N\left(
f_{i_{1}},f_{i_{2}},...,f_{i_{l}},0,0,...,0\right) =$

\ 

$\sum\limits_{\left( X_{i}\right) _{i\in I_{\nu ^{\prime }}}\in \mathbb{F}%
_{q}^{\left| I_{\nu ^{\prime }}\right| }}\prod\limits_{t=1}^{l}\left(
q^{\left| J_{\nu _{t}}\right| -1}+\kappa \left(
\sum\limits_{j=1}^{m-t}\sigma _{\left[ i_{1},i_{2},...,i_{l}\right]
}^{t,j}\prod\limits_{\tau \in J_{\nu _{j}^{,}}}X_{\tau }\right) q^{\left|
J_{\nu _{t}}\right| -1}\left( 1-\left( \frac{q-1}{q}\right) ^{\left| J_{\nu
_{t}}\right| -1}\right) \right) .$

\ 

After, we get the desired result for the case $a=0$ from $\left( 2.6b\right)
.$

\ 

Now assume that $a\neq 0$ and let \ $\left( a_{1},...,a_{k}\right) \in 
\mathbb{F}_{q}^{k}$ such that $a_{1}a_{2}...a_{k}=a,$ then from
Proposition~III-3, we can write:

\ 

$N\left( f_{1},f_{2},...,f_{k},a_{1},....,a_{k}\right) =$

\ 

$\sum\limits_{\left( X_{i}\right) _{i\in I_{\nu ^{\prime }}}\in \mathbb{F}%
_{q}^{\left| I_{\nu ^{\prime }}\right| }}\prod\limits_{t=1}^{k}\left(
q^{\left| J_{\nu _{t}}\right| -1}+\kappa \left(
\sum\limits_{j=1}^{k}a_{t,j}^{^{\prime }}a_{j}-\sum\limits_{j=1}^{m-k}\sigma
_{\left[ 1,2,...,k\right] }^{t,j}\prod\limits_{\tau \in J_{\nu
_{j}^{,}}}X_{\tau }\right) q^{\left| J_{\nu _{i}}\right| -1}\left( 1-\left( 
\frac{q-1}{q}\right) ^{\left| J_{\nu _{i}}\right| -1}\right) \right) ,$

\ 

it follows from $\left( 2.6a\right) $ that:

\ 

$N\left( f,a\right) =$

$\sum\limits_{\substack{ \left( a_{1},...,a_{k}\right) \in \mathbb{F}%
_{q}^{k}  \\ a_{1}...a_{k}=a}}\sum\limits_{\left( X_{i}\right) _{i\in I_{\nu
^{\prime }}}\in \mathbb{F}_{q}^{\left| I_{\nu ^{\prime }}\right|
}}\prod\limits_{t=1}^{k}(q^{\left| J_{\nu _{t}}\right| -1}\left( \text{1+}%
\kappa \left( \sum\limits_{j=1}^{k}a_{t,j}^{^{\prime
}}a_{j}-\sum\limits_{j=1}^{m-k}\sigma _{\left[ 1,2,...,k\right]
}^{t,j}\prod\limits_{\tau \in J_{\nu _{j}^{,}}}X_{\tau }\right) \left( \text{%
1-}\left( \frac{q-1}{q}\right) ^{\left| J_{\nu _{t}}\right| -1}\right)
\right) ),$

\ 

which completes the proof.$\blacksquare $

\ 

The formula given in previous Proposition is generally a semi-explicit
formula because of\ terms 
\begin{equation*}
\kappa \left( \sum\limits_{j=1}^{m-t}\sigma _{\left[ i_{1},i_{2},...,i_{l}%
\right] }^{t,j}\prod\limits_{\tau \in J_{\nu _{j}^{,}}}X_{\tau }\right) ,\ \
\kappa \left( \sum\limits_{j=1}^{k}a_{t,j}^{^{\prime
}}a_{j}-\sum\limits_{j=1}^{m-k}\sigma _{\left[ 1,2,...,k\right]
}^{t,j}\prod\limits_{\tau \in J_{\nu _{j}^{,}}}X_{\tau }\right) .
\end{equation*}
In the next Proposition, particular cases related to the matrix $A$ are
considered.

\ 

\textbf{PROPOSITION IV-2}

$\,$

$If$ $m=k$ $and$ $a\neq 0$ $then$:

\ 

$\left( 4.1\right) $\ \ \ \ \ \ \ \ \ \ \ $N\left( f,a\right) =$\ 

\ 

$\sum\limits_{\left( a_{1},..,a_{m-1}\right) \in \left( \mathbb{F}_{q}^{\ast
}\right) ^{m-1}}[\prod\limits_{i=1}^{m}\left( q^{\left| J_{i}\right|
-1}+\kappa \left( \sum\limits_{j=1}^{m-1}a_{ij}^{^{\prime
}}a_{j}+a_{im}^{^{\prime }}\frac{a}{a_{1}....a_{m-1}}\right) q^{\left|
J_{i}\right| -1}\left( 1-\left( \frac{q-1}{q}\right) ^{\left| J_{i}\right|
-1}\right) \right) ].$

\ 

$If$ \ 

$\left( 4.2\right) \ \ \ \ \ \ \ \ \ \ \ \ \ \ \ \ \ \ \ \ \ \ \left\{ 
\begin{array}{c}
m=2k, \\ 
A=\left( D_{1}\text{ }D_{2}\right) ,%
\end{array}
\right. $

\ 

$where$ $D_{1}$ $and$ $D_{2}$\ $are$ $two$ $invertible$ $diagonale$ $k\times
k$ $matrices,$ $then$:

\ \ 

$\left( 4.3\right) \ \ \ \ N\left( f,a\right) =$\ 

\ 

\ \ \ \ \ \ \ \ $\left( q-1\right) ^{k-1}$\ $\prod\limits_{j=1}^{k}\left[
q^{\left| J_{i}\right| +\left| J_{k+i}\right| -1}-q^{\left| J_{i}\right|
+\left| J_{k+i}\right| -1}\left( 1-\left( \frac{q-1}{q}\right) ^{\left|
J_{i}\right| -1}\right) \left( 1-\left( \frac{q-1}{q}\right) ^{\left|
J_{k+i}\right| -1}\right) \right] ,$

\ \ \ \ \ \ \ \ \ \ $\ if$ \ $a\neq 0,$

\ 

$\left( 4.4\right) $ $\ \ \ N\left( f,0\right) $\ $=$\ $kq^{n-1}+\left(
q-1\right) q^{n-1}\sum\limits_{i=1}^{k}[\left( 1-\left( \frac{q-1}{q}\right)
^{\left| J_{i}\right| -1}\right) \left( 1-\left( \frac{q-1}{q}\right)
^{\left| J_{k+i}\right| -1}\right) ]-$

\ 

$\sum\limits_{1\leq i_{1}<i_{2}\leq k}q^{n-\sum\limits_{j=1}^{2}(\left|
J_{i_{j}}\right| +\left| J_{k+i_{j}}\right|
)}\prod\limits_{j=1}^{2}(q^{\left| J_{i_{j}}\right| +\left|
J_{k+i_{j}}\right| -1}\left( \text{1+}\kappa \left( 0\right) \left( 1-\left( 
\frac{q-1}{q}\right) ^{\left| J_{i_{j}}\right| -1}\right) \left( 1-\left( 
\frac{q-1}{q}\right) ^{\left| J_{k+i_{j}}\right| -1}\right) \right) )+$

$...+$

$\left( -1\right) ^{l+1}\sum\limits_{1\leq i_{1}<..<i_{l}\leq
k}q^{n-\sum\limits_{j=1}^{l}(\left| J_{i_{j}}\right| +\left|
J_{k+i_{j}}\right| )}\prod\limits_{j=1}^{l}(q^{\left| J_{i_{j}}\right|
+\left| J_{k+i_{j}}\right| -1}\left( \text{1+}\kappa \left( 0\right) \left( 
\text{1-}\left( \frac{q-1}{q}\right) ^{\left| J_{i_{j}}\right| -1}\right)
\left( 1\text{-}\left( \frac{q-1}{q}\right) ^{\left| J_{k+i_{j}}\right|
-1}\right) \right) )$

\ 

$+...+\left( -1\right) ^{k+1}\prod\limits_{j=1}^{k}(q^{\left|
J_{i_{j}}\right| +\left| J_{k+i_{j}}\right| -1}\left( 1+\kappa \left(
0\right) \left( 1-\left( \frac{q-1}{q}\right) ^{\left| J_{i_{j}}\right|
-1}\right) \left( 1-\left( \frac{q-1}{q}\right) ^{\left| J_{k+i_{j}}\right|
-1}\right) \right) ).$

\ \ 

\textbf{Proof}

\ 

First consider the case $m=k.$ Then $A$ is an invertible $m\times m$ matrix.

Let\ $\left( a_{1},a_{2},...,a_{m-1}\right) \in \left( \mathbb{F}_{q}^{\ast
}\right) ^{m-1}.$ In\ $\left( 2.4a\right) ,$ $\left( 2.5a\right) $ we can
take:

\ 

\ \ \ \ \ \ \ \ \ $B_{\left[ i_{1},i_{2},...,i_{l}\right] }=A,$ $C_{\left[
i_{1},i_{2},...,i_{l}\right] }=0,$ $\left( b_{1},...,b_{m}\right) =\left(
a_{1},a_{2},...,a_{m-1},\frac{a}{a_{1}a_{2}....a_{m-1}}\right) ,$

\ 

\ \ \ \ \ \ \ \ \ $l=m,$ $\left[ i_{1},i_{2},...,i_{m}\right] =\left[
1,2,..,m\right] ,$ $I_{\nu ^{^{\prime }}}=\emptyset ,$

\ 

then from $\left( 3.2\right) $, we can write:

\ 

$N\left( f_{1},f_{2},...,f_{m},a_{1},....,a_{m-1},\frac{a}{a_{1}....a_{m-1}}%
\right) =$

\ \ \ \ \ \ \ \ 

\ \ \ \ \ \ \ \ \ \ \ \ \ \ \ \ \ \ \ \ \ \ \ \ $\prod\limits_{i=1}^{m}%
\left( q^{\left| J_{i}\right| -1}+\kappa \left(
\sum\limits_{j=1}^{m-1}a_{i,j}^{^{\prime }}a_{j}+a_{im}^{^{\prime }}\frac{a}{%
a_{1}...a_{m-1}}\right) q^{\left| J_{i}\right| -1}\left( 1-\left( \frac{q-1}{%
q}\right) ^{\left| J_{i}\right| -1}\right) \right) $

\ 

and the result $\left( 4.1\right) $ follows from $\left( 2.6a\right) $.

\ 

Now consider the case $\left( 4.2\right) \ $and let $\left(
a_{1},a_{2},...,a_{k-1}\right) \in \left( \mathbb{F}_{q}^{\ast }\right)
^{k-1}.$ Recalling that $A=\left( a_{ij}\right) ,$ it's easy to see that:\ $%
D_{1}=Diag\left( a_{11},...,a_{kk}\right) $, $D_{2}=Diag\left(
a_{1,k+1},a_{2,k+2}...,a_{k,k+k}\right) $, $a_{ij}=0$ \ for all $\left(
i,j\right) $\ such\ that $j-i\notin \left\{ 0,i\right\} ,\
a_{ii}a_{i,k+i}\neq 0$\ for all $i\in \left\{ 1,...,k\right\} $. This means
that we can make the following choices in $\left( 2.4a\right) $:

\ 

\ \ \ \ \ \ \ \ \ $B_{\left[ i_{1},i_{2},...,i_{l}\right] }=D_{1},$ $C_{%
\left[ i_{1},i_{2},...,i_{l}\right] }=D_{2},$ $\left( b_{1},...,b_{k}\right)
^{T}=\left( a_{1},a_{2},...,a_{k-1},\frac{a}{a_{1}a_{2}....a_{k-1}}\right)
^{T},$

\ 

\ \ \ \ \ \ \ \ \ $l=k,$ $\left[ i_{1},i_{2},...,i_{k}\right] =\left[
1,2,..,k\right] ,$ $\left( \nu _{1},..,\nu _{k}\right) =\left(
1,2,..,k\right) ,$

\ 

$\ \ \ \ \ \ \ \ \ \left( \nu _{1}^{^{\prime }},..,\nu _{m-k}^{^{\prime
}}\right) =\left( k+1,..,m-k\right) ,$

\ \ \ 

thus:\ \ \ \ \ 

\ \ \ \ 

$\ \ \ \ \left( \prod\limits_{\tau \in J_{1}}X_{\tau },\prod\limits_{\tau
\in J_{2}}X_{\tau },....,\prod\limits_{\tau \in J_{k}}X_{\tau }\right)
^{T}=D_{1}^{-1}\left( a_{1},...,a_{k-1},\frac{a}{a_{1}a_{2}....a_{k-1}}%
\right) ^{T}-$

\ \ 

$\ \ \ \ \ \ \ \ \ \ \ \ \ \ \ \ \ \ \ \ \ \ \ \ \ \ \ \ \ \ \ \ \ \ \ \ \ \
\ \ \ \ \ \ \ \ \ \ \ \ \ \ \ \ \ D_{1}^{-1}D_{2}\left( \prod\limits_{\tau
\in J_{k+1}}X_{\tau },\prod\limits_{\tau \in J_{k+2}}X_{\tau
},....,\prod\limits_{\tau \in J_{m}}X_{\tau }\right) ^{T}$

\ 

and for all $i\in \left\{ 1,..,k\right\} $,

\ \ 

$\left( 4.5\right) $\ \ \ \ \ \ \ \ \ \ \ \ \ \ $\prod\limits_{\tau \in
J_{i}}X_{\tau }=a_{ii}^{-1}a_{i}-a_{ii}^{-1}a_{i_{,}k+i}\prod\limits_{\tau
\in J_{k+i}}X_{\tau },$

\ \ \ \ \ \ \ \ \ \ \ \ \ \ \ \ \ \ with \ $a_{k}=\frac{a}{a_{1}....a_{k-1}}%
. $

\ \ 

We deduce from Proposition III-1 and the hypothesis $\left(
a_{ii},a_{i},a_{i_{,}k+i}\right) \in \left( \mathbb{F}_{q}^{\ast }\right)
^{3}$, that\ the\ number of solutions in $\mathbb{F}_{q}^{\left|
J_{i}\right| +\left| J_{k+i}\right| }$ of $\left( 4.5\right) $ is:

\ 

$\ \ \ \ \ \ \ \ \ \ \ \ \ q^{\left| J_{i}\right| +\left| J_{k+i}\right|
-1}-q^{\left| J_{i}\right| +\left| J_{k+i}\right| -1}\left( 1-\left( \frac{%
q-1}{q}\right) ^{\left| J_{i}\right| -1}\right) \left( 1-\left( \frac{q-1}{q}%
\right) ^{\left| J_{k+i}\right| -1}\right) .$

\ 

Then, it's easy to see that:

\ 

\ \ \ \ $N\left( f_{1},f_{2},...,f_{k},a_{1},....,a_{k-1},\frac{a}{%
a_{1}a_{2}....a_{k-1}}\right) =$

\ 

\ \ \ \ \ \ \ \ \ \ \ \ \ $\prod\limits_{i=1}^{k}\left[ q^{\left|
J_{i}\right| +\left| J_{k+i}\right| -1}-q^{\left| J_{i}\right| +\left|
J_{k+i}\right| -1}\left( 1-\left( \frac{q-1}{q}\right) ^{\left| J_{i}\right|
-1}\right) \left( 1-\left( \frac{q-1}{q}\right) ^{\left| J_{k+i}\right|
-1}\right) \right] $

\ 

and from $\left( 2.6a\right) $ the result $\left( 4.3\right) $ follows.

\ \ 

Now, let $1\leq i_{1}<i_{2}<...<i_{l}\leq k.$ The system of equations $%
\left( 2.3\right) $ in the case $\left( 4.2\right) $ and $\left(
b_{1},...,b_{l}\right) =\left( 0,...,0\right) $, namely:

\ \ 

\ \ \ \ \ \ \ \ \ \ \ \ \ \ \ \ \ \ \ \ \ \ \ \ \ \ $\left\{ 
\begin{array}{c}
\ f_{i_{1}}\left( X_{1},..,X_{n}\right) =0, \\ 
\ f_{i_{2}}\left( X_{1},..,X_{n}\right) =0, \\ 
... \\ 
f_{i_{l}}\left( X_{1},..,X_{n}\right) =0,%
\end{array}
\right. $

\ 

is equivalent to

\ \ \ \ \ \ \ \ \ \ \ \ \ \ \ \ \ \ \ \ \ \ \ \ \ \ $\left\{ 
\begin{array}{c}
\ \ \prod\limits_{\tau \in J_{i_{1}}}X_{\tau
}+a_{i_{1}i_{1}}^{-1}a_{i_{1,}k+i_{1}}\prod\limits_{\tau \in
J_{k+i_{1}}}X_{\tau }=0, \\ 
\prod\limits_{\tau \in J_{i_{2}}}X_{\tau
}+a_{i_{2}i_{2}}^{-1}a_{i_{2,}k+i_{2}}\prod\limits_{\tau \in
J_{k+i_{2}}}X_{\tau }=0, \\ 
... \\ 
\prod\limits_{\tau \in J_{i_{l}}}X_{\tau
}+a_{i_{l}i_{l}}^{-1}a_{i_{l,}k+i_{l}}\prod\limits_{\tau \in
J_{k+i_{l}}}X_{\tau }=0,%
\end{array}
\right. $

thus:

\ 

\ \ \ \ \ \ \ \ \ \ \ \ \ \ $N\left( \
f_{i_{1}},f_{i_{2}},...,f_{i_{l}},0,0,..,0\right) =$

\ 

$q^{n-\sum\limits_{j=1}^{l}(\left| J_{i_{j}}\right| +\left|
J_{k+i_{j}}\right| )}\prod\limits_{j=1}^{l}\left( q^{\left| J_{i_{j}}\right|
+\left| J_{k+i_{j}}\right| -1}+\kappa \left( 0\right) q^{\left|
J_{i_{j}}\right| +\left| J_{k+i_{j}}\right| -1}\left( 1\text{-}\left( \frac{%
q-1}{q}\right) ^{\left| J_{i_{j}}\right| -1}\right) \left( 1-\left( \frac{q-1%
}{q}\right) ^{\left| J_{k+i_{j}}\right| -1}\right) \right) .$

After by Proposition III-1$,$ we have for all \ $j\in \left\{ 1,..,l\right\} 
$:

\ 

\ \ \ \ \ \ $N\left( f_{i},0\right) $\ $=q^{n-1}+\kappa \left( 0\right)
q^{n-1}\left( 1-\left( \frac{q-1}{q}\right) ^{\left| J_{i}\right| -1}\right)
\left( 1-\left( \frac{q-1}{q}\right) ^{\left| J_{k+i}\right| -1}\right) ,$

\ \ 

consequently from $\left( 2.6b\right) ,$ we have:

\ 

$N\left( f,0\right) $\ $=$\ $kq^{n-1}+\kappa \left( 0\right)
q^{n-1}\sum\limits_{i=1}^{k}[\left( 1-\left( \frac{q-1}{q}\right) ^{\left|
J_{i}\right| -1}\right) \left( 1-\left( \frac{q-1}{q}\right) ^{\left|
J_{k+i}\right| -1}\right) ]-$

\ 

$\sum\limits_{1\leq i_{1}<i_{2}\leq k}q^{n-\sum\limits_{j=1}^{2}(\left|
J_{i_{j}}\right| +\left| J_{k+i_{j}}\right| )}\left[ \prod%
\limits_{j=1}^{2}(q^{\left| J_{i_{j}}\right| +\left| J_{k+i_{j}}\right| 
\text{-}1}\left( \text{1+}\kappa \left( 0\right) \left( \text{1-}\left( 
\frac{q-1}{q}\right) ^{\left| J_{i_{j}}\right| \text{-}1}\right) \left( 
\text{1-}\left( \frac{q-1}{q}\right) ^{\left| J_{k+i_{j}}\right| \text{-}%
1}\right) \right) )\right] +...+$

\ 

$\left( \text{-}1\right) ^{l+1}\sum\limits_{1\leq i_{1}<...<i_{l}\leq k}q^{n%
\text{-}\sum\limits_{j=1}^{l}(\left| J_{i_{j}}\right| \text{+}\left|
J_{k+i_{j}}\right| )}\prod\limits_{j=1}^{l}(q^{\left| J_{i_{j}}\right| \text{%
+}\left| J_{k+i_{j}}\right| \text{-}1}\left( \text{1+}\kappa \left( 0\right)
\left( \text{1-}\left( \frac{q\text{-}1}{q}\right) ^{\left| J_{i_{j}}\right| 
\text{-}1}\right) \left( \text{1-}\left( \frac{q\text{-}1}{q}\right)
^{\left| J_{k+i_{j}}\right| \text{-}1}\right) \right) )$

\ 

$+...+\left( -1\right) ^{k+1}\prod\limits_{j=1}^{k}(q^{\left|
J_{i_{j}}\right| +\left| J_{k+i_{j}}\right| -1}\left( \text{1+}\kappa \left(
0\right) \left( \text{1}-\left( \frac{q-1}{q}\right) ^{\left|
J_{i_{j}}\right| -1}\right) \left( \text{1}-\left( \frac{q-1}{q}\right)
^{\left| J_{k+i_{j}}\right| -1}\right) \right) ),$

\ \ 

which is the formula $\left( 4.4\right) .\blacksquare $

\ 

This Proposition indicates for instance that for particular interesting
cases related to the matrix~$A$, we can get an explicit formula for the
number of solutions of $\left( 1.3\right) .$\bigskip \bigskip

V \ NUMERICAL EXAMPLES

\ \ 

\ \ 

\ \ V-1 \ A FIRST EXAMPLE

\ \ 

Consider the following particular situation of $\left( 1.3\right) $:

\ 

$\left( 5.1\right) $\ \ $\ \ \ \ \ \ \ \left(
X_{1}X_{2}+X_{5}X_{6}X_{7}\right) \left( X_{3}X_{4}+X_{5}X_{6}X_{7}\right) $%
\ $=a$\ .

\ 

\ \ Then:\ \ \ \ \ \ \ \ \ \ \ 

\ \ 

\ $\ \ \ \left\{ 
\begin{array}{c}
\left( k,m,n\right) =\left( 2,3,7\right) ,\smallskip \\ 
A=\left( 
\begin{array}{ccc}
1 & 0 & 1 \\ 
0 & 1 & 1%
\end{array}
\right) ,\smallskip \\ 
\ J_{1}=\left\{ 1,2\right\} ,\text{ }J_{2}=\left\{ 3,4\right\} ,\
J_{3}=\left\{ 5,6,7\right\} .%
\end{array}
\right. $

\ \ \ \ 

\ We have the following Proposition.

\ 

\textbf{PROPOSITION V-1}

$\;$

$Let$ $\alpha $ $be$ $a$ $primitive$ $element$ $of$ $\mathbb{F}_{q},$ $1\leq
\tau \leq \frac{q-1}{2}.$ $Then$ $the$ $number$ $of$ $solutions$\ \ $of$ \ $%
\left( 5.1\right) $ $is$:\ 

\ \ 

$\left( 5.2\right) $\ \ \ \ $N\left( f,0\right)
=2q^{6}+3q^{5}-13q^{4}+16q^{3}-9q^{2}+2q,$

\ \ \ 

$\left( 5.3\right) $\ \ $\ \ N\left( f,\alpha ^{2\tau }\right) =\left(
q-1\right) ^{2}q\left( q^{3}+q^{2}-2q+2\right) ,$ \ \ \ \ 

\ 

$\left( 5.4\right) $ $\ \ \ N\left( f,\alpha ^{2\tau -1}\right) =\left(
q-1\right) ^{3}q\left( q^{2}+2q-2\right) .$

\thinspace

\textbf{Proof}

\thinspace

In $\ \left( 2.4a\right) $ we can take:

$\ l=2,$ $\left[ i_{1},i_{2}\right] =\left[ 1,2\right] ,$ $\left( \nu
_{1},\nu _{2}\right) =\left( 1,2\right) ,$ $\nu _{1}^{^{\prime }}=3,$\ $B_{%
\left[ 1,2\right] }=\left( 
\begin{array}{cc}
1 & 0 \\ 
0 & 1%
\end{array}
\right) ,$ $C_{\left[ 1,2\right] }=\left( 
\begin{array}{c}
1 \\ 
1%
\end{array}
\right) ,$

thus:

\ 

$\ B_{\left[ 1,2\right] }^{-1}=\left( 
\begin{array}{cc}
1 & 0 \\ 
0 & 1%
\end{array}
\right) =\left( 
\begin{array}{cc}
a_{11}^{^{\prime }} & a_{12}^{^{\prime }} \\ 
a_{21}^{^{\prime }} & a_{22}^{^{\prime }}%
\end{array}
\right) ,$ $B_{\left[ 1,2\right] }^{-1}$ $C_{\left[ 1,2\right] }=\left( 
\begin{array}{c}
1 \\ 
1%
\end{array}
\right) =\left( 
\begin{array}{c}
\sigma _{\left[ 1,2\right] }^{1,1} \\ 
\sigma _{\left[ 1,2\right] }^{2,1}%
\end{array}
\right) .$

\ 

From the above assumptions and Proposition III-4, we first can establish
easily $\left( 5.2\right) $ and secondly we can deduce that for all $\
1<t\leq q-1$:

\ \ \ \ \ \ \ \ 

$\left( 5.5\right) $ $N\left( f,\text{ }\alpha ^{t}\right)
=\sum\limits_{u\in \mathbb{F}_{q}^{\ast }}\sum\limits_{\left(
X_{5},X_{6},X_{7}\right) \in \mathbb{F}_{q}^{3}}\left[ \left( q+\kappa
\left( u-X_{5}X_{6}X_{7}\right) \right) \left( q+\kappa \left( \frac{\alpha
^{t}}{u}-X_{5}X_{6}X_{7}\right) \right) \right] ,$

\ 

it follows that:\ \ \ 

\ 

$N\left( f,\text{ }\alpha ^{2\tau }\right) =\sum\limits_{u\in \left\{
-\alpha ^{\tau },\alpha ^{\tau }\right\} }\sum\limits_{\left(
X_{5},X_{6},X_{7}\right) \in \mathbb{F}_{q}^{3}}\left[ \left( q+\kappa
\left( u-X_{5}X_{6}X_{7}\right) \right) \left( q+\kappa \left( \frac{\alpha
^{2\tau }}{u}-X_{5}X_{6}X_{7}\right) \right) \right] +$

\ 

\ \ \ \ \ \ \ \ \ \ \ \ \ \ \ \ \ \ \ $\sum\limits_{u\notin \left\{ -\alpha
^{\tau },0,\text{ }\alpha ^{\tau }\right\} }\sum\limits_{\left(
X_{5},X_{6},X_{7}\right) \in \mathbb{F}_{q}^{3}}\left[ \left( q+\kappa
\left( u-X_{5}X_{6}X_{7}\right) \right) \left( q+\kappa \left( \frac{\alpha
^{2\tau }}{u}-X_{5}X_{6}X_{7}\right) \right) \right] $

\ 

\ \ \ \ \ \ \ \ \ \ \ \ \ \ \ \ $=T_{1}+T_{2},$

\ 

where:

\ 

\ $\left( 5.6\right) $\ \ \ \ \ \ \ \ $T_{1}=\sum\limits_{u\in \left\{
-\alpha ^{\tau },\alpha ^{\tau }\right\} }\sum\limits_{\left(
X_{5},X_{6},X_{7}\right) \in \mathbb{F}_{q}^{3}}\left[ \left( q+\kappa
\left( u-X_{5}X_{6}X_{7}\right) \right) \left( q+\kappa \left( \frac{\alpha
^{2\tau }}{u}-X_{5}X_{6}X_{7}\right) \right) \right] ,$

\ $\left( 5.7\right) \ \ \ \ \ \ \ \ T_{2}=$\ \ $\sum\limits_{u\notin
\left\{ -\alpha ^{\tau },0,\text{ }\alpha ^{\tau }\right\}
}\sum\limits_{\left( X_{5},X_{6},X_{7}\right) \in \mathbb{F}_{q}^{3}}\left[
\left( q+\kappa \left( u-X_{5}X_{6}X_{7}\right) \right) \left( q+\kappa
\left( \frac{\alpha ^{2\tau }}{u}-X_{5}X_{6}X_{7}\right) \right) \right] ,$

\ 

thus:

\ 

\ \ \ $T_{1}=$\ \ $\sum\limits_{\left( X_{5},X_{6},X_{7}\right) \in \mathbb{F%
}_{q}^{3}}\left( q+\kappa \left( -\alpha ^{\tau }-X_{5}X_{6}X_{7}\right)
\right) ^{2}$\ +\ \ $\sum\limits_{\left( X_{5},X_{6},X_{7}\right) \in 
\mathbb{F}_{q}^{3}}\left( q+\kappa \left( \alpha ^{\tau
}-X_{5}X_{6}X_{7}\right) \right) ^{2}$\ \ 

\ 

\ \ \ \ \ \ $=\sum\limits_{\substack{ \left( X_{5},X_{6},X_{7}\right) \in 
\mathbb{F}_{q}^{3}/  \\ X_{5}X_{6}X_{7}=-\alpha ^{\tau }}}\left( q+\kappa
\left( 0\right) \right) ^{2}$\ \ +\ $\sum\limits_{\substack{ \left(
X_{5},X_{6},X_{7}\right) \in \mathbb{F}_{q}^{3}/  \\ X_{5}X_{6}X_{7}\neq
-\alpha ^{\tau }}}\left( q+\kappa \left( -\alpha ^{\tau
}-X_{5}X_{6}X_{7}\right) \right) ^{2}$

\ 

\ \ \ \ \ \ \ \ \ \ \ \ \ \ \ +\ $\sum\limits_{\substack{ \left(
X_{5},X_{6},X_{7}\right) \in \mathbb{F}_{q}^{3}/  \\ X_{5}X_{6}X_{7}=\alpha
^{\tau }}}\left( q+\kappa \left( 0\right) \right) ^{2}+\sum\limits 
_{\substack{ \left( X_{5},X_{6},X_{7}\right) \in \mathbb{F}_{q}^{3}/  \\ %
X_{5}X_{6}X_{7}\neq \alpha ^{\tau }}}\left( q+\kappa \left( \alpha ^{\tau
}-X_{5}X_{6}X_{7}\right) \right) ^{2}.$

\ 

Using Corollary \ III-2, we get:

\ 

$\left( 5.8\right) $ \ \ \ \ $T_{1}=2\left( q-1\right) ^{2}\left(
2q-1\right) ^{2}+2\left( q^{3}-\left( q-1\right) ^{2}\right) \left(
q-1\right) ^{2}$

\ \ \ \ \ \ \ \ \ \ \ \ \ \ \ \ \ \ \ $=2\left( q-1\right) ^{2}q\left[
q^{2}+3q-2\right] $.

Now,\ $\left( 5.7\right) $ implies that:

\ 

\ \ \ \ \ $T_{2}=\sum\limits_{u\notin \left\{ -\alpha ^{\tau },0,\text{ }%
\alpha ^{\tau }\right\} }[\sum\limits_{\substack{ \left(
X_{5},X_{6},X_{7}\right) \in \mathbb{F}_{q}^{3}/  \\ X_{5}X_{6}X_{7}=u}}%
\left[ \left( q+\kappa (0\right) )\left( q+\kappa \left( \frac{\alpha
^{2\tau }}{u}-X_{5}X_{6}X_{7}\right) \right) \right] +$

\ 

\ $\sum\limits_{\substack{ \left( X_{5},X_{6},X_{7}\right) \in \mathbb{F}%
_{q}^{3}/  \\ X_{5}X_{6}X_{7}=\frac{\alpha ^{2\tau }}{u}}}\left[ \left(
q+\kappa \left( u-X_{5}X_{6}X_{7}\right) \right) \left( q+\kappa \left(
0\right) \right) \right] $\ $+$

\ \ \ 

\ \ \ \ \ \ $\sum\limits_{\substack{ \left( X_{5},X_{6},X_{7}\right) \in 
\mathbb{F}_{q}^{3}/  \\ X_{5}X_{6}X_{7}\notin \left\{ u,\frac{\alpha ^{2\tau
}}{u}\right\} }}\left[ \left( q+\kappa \left( u-X_{5}X_{6}X_{7}\right)
\right) \left( q+\kappa \left( \frac{\alpha ^{2\tau }}{u}-X_{5}X_{6}X_{7}%
\right) \right) \right] ].$

\ 

It follows from Corollary III-2 that:

\ 

\ \ \ $\left( 5.9\right) $\ \ \ $T_{2}=\sum\limits_{u\notin \left\{ -\alpha
^{\tau },0,\text{ }\alpha ^{\tau }\right\} }\left[ 2\left( q-1\right)
^{3}\left( 2q-1\right) +\left( q^{3}-2\left( q-1\right) ^{2}\right) \left(
q-1\right) ^{2}\right] $

\ \ \ \ \ \ \ \ \ \ \ \ \ \ \ \ \ \ \ $=\left( q-3\right) \left[ 2\left(
q-1\right) ^{3}\left( 2q-1\right) +\left( q^{3}-2\left( q-1\right)
^{2}\right) \left( q-1\right) ^{2}\right] $

\ 

\ \ \ \ \ \ \ \ \ \ \ \ \ \ \ \ \ \ \ $=q\left( q-3\right) \left( q-1\right)
^{2}\left( q^{2}+2q-2\right) ,$

\ \ \ 

and from $\left( 5.8\right) $ and $\left( 5.9\right) $, the formula $\left(
5.3\right) $ follows.

\ 

Now a similar analysis leads to the following:

\ 

$N\left( f,\text{ }\alpha ^{2\tau -1}\right) =\sum\limits_{u\in \mathbb{F}%
_{q}^{\ast }}[\sum\limits_{\substack{ \left( X_{5},X_{6},X_{7}\right) \in 
\mathbb{F}_{q}^{3}/  \\ X_{5}X_{6}X_{7}=u}}\left[ \left( q+\kappa (0\right)
)\left( q+\kappa \left( \frac{\alpha ^{2\tau -1}}{u}-X_{5}X_{6}X_{7}\right)
\right) \right] +$

\ 

\ $\sum\limits_{\substack{ \left( X_{5},X_{6},X_{7}\right) \in \mathbb{F}%
_{q}^{3}/  \\ X_{5}X_{6}X_{7}=\frac{\alpha ^{2\tau -1}}{u}}}\left[ \left(
q+\kappa \left( u-X_{5}X_{6}X_{7}\right) \right) \left( q+\kappa \left(
0\right) \right) \right] $\ $+$

\ \ \ 

\ \ \ \ \ \ $\sum\limits_{\substack{ \left( X_{5},X_{6},X_{7}\right) \in 
\mathbb{F}_{q}^{3}/  \\ X_{5}X_{6}X_{7}\notin \left\{ u,\frac{\alpha ^{2\tau
-1}}{u}\right\} }}\left[ \left( q+\kappa \left( u-X_{5}X_{6}X_{7}\right)
\right) \left( q+\kappa \left( \frac{\alpha ^{2\tau -1}}{u}%
-X_{5}X_{6}X_{7}\right) \right) \right] ,$

\ 

\ \ \ \ \ \ \ \ \ \ \ \ \ \ \ \ $=\left( q-1\right) \left[ 2\left(
q-1\right) ^{3}\left( 2q-1\right) +\left( q^{3}-2\left( q-1\right)
^{2}\right) \left( q-1\right) ^{2}\right] .$

\ 

It's follows that $N\left( f,\text{ }\alpha ^{2\tau -1}\right) =\left(
q-1\right) ^{3}q\left( q^{2}+2q-2\right) $, which is the formula $\left(
5.4\right) .\blacksquare $

\ 

The result given in this Proposition indicates for particular cases of $%
\left( 1.3\right) $ that the obtained solution is a simple explicit formula.

\ \ 

\ \ 

V-2 \ A SECOND EXAMPLE

\ \ 

The example that we consider here shows us that the obtained formula for $%
N\left( f,a\right) $ in Proposition IV-1 is more efficient than a trivial
systematic computation of all the solutions. In the semi-explicit formula $%
\left( 5.12\right) $, the use of a classical machine is certainly required,
but lead to easier computation involving sums and the obvious function $%
\kappa $. This example is the following:

\ 

$\left( 5.10\right) $ \ \ \ $\left( X_{1}+X_{6}X_{7}X_{8}\right) \left(
X_{1}+X_{2}X_{3}+X_{6}X_{7}X_{8}\right) \left( X_{2}X_{3}+X_{4}X_{5}\right)
=a$

\ 

We have:

\ \ 

$\ \left\{ 
\begin{array}{c}
\left( k,m,n\right) =\left( 3,4,8\right) ,\smallskip \\ 
A=\left( 
\begin{array}{cccc}
1 & 0 & 0 & 1 \\ 
1 & 1 & 0 & 1 \\ 
0 & 1 & 1 & 0%
\end{array}
\right) ,\smallskip \\ 
J_{1}=\left\{ 1\right\} ,\text{ }J_{2}=\left\{ 2,3\right\} ,\ J_{3}=\left\{
4,5\right\} ,\text{ }J_{4}=\left\{ 6,7,8\right\}%
\end{array}
\right. $

\ 

and the following result:

\bigskip

\textbf{PROPOSITION V-2}

\ 

$The$ $number$ $N\left( f,a\right) $ $of$ $solutions$\ \ $of$ \ $\left(
5.10\right) $ $is$:

\ 

$\left( 5.11\right) $ \ \ \ \ $N\left( f,0\right) =3q^{7}-3q^{6}+q^{3}\left(
2q^{2}-2q+1\right) ,$

\ 

$\left( 5.12\right) $ \ \ \ \ $N\left( f,a\right)
=q^{3}\sum\limits_{a_{1},a_{2}\in \mathbb{F}_{q}^{\ast }}\left( q+\kappa
\left( a_{2}-a_{1}\right) \right) \left( q+\kappa \left( a_{1}-a_{2}+\frac{a%
}{a_{1}a_{2}}\right) \right) ,$

\ \ \ \ \ \ \ \ \ \ \ \ \ \ \ $if$ $a\neq 0.$

\ 

\textbf{Proof}

\thinspace

With\ obvious notations, Proposition IV-1 gives:

\ 

$\left( 5.13\right) $\ \ $\ \ \ N\left( f,0\right) =\xi -\
\sum\limits_{1\leq i_{1}<i_{2}\leq 3}$\ $C_{i_{1},i_{2}}+C_{1,2,3}.%
\smallskip $

It is easily seen that $\xi =3q^{7}+\left( q-1\right) q^{5}$. To compute $%
C_{1,2},$ we consider in $\left( 2.4a\right) $:

\ \ 

\ \ \ \ $B_{\left[ 1,2\right] }=\left( 
\begin{array}{cc}
1 & 0 \\ 
1 & 1%
\end{array}
\right) $ , $\ C_{\left[ 1,2\right] }=\left( 
\begin{array}{cc}
0 & 1 \\ 
0 & 1%
\end{array}
\right) .$

\ 

Hence $B_{\left[ 1,2\right] }^{-1}$ $C_{\left[ 1,2\right] }=\left( 
\begin{array}{cc}
0 & 1 \\ 
0 & 0%
\end{array}
\right) =\left( \sigma _{\left[ 1,2\right] }^{i,j}\right) _{i,j}$ and \ $%
C_{1,2}=q^{5}\left( 2q-1\right) ,$ as is easy to check. Likewise we obtain:

\ \ 

$C_{1,3}=C_{2,3}=q^{6}+q^{3}\sum\limits_{\left( X_{4},X_{5}\right) \in
F_{q}^{2}}\kappa \left( X_{4}X_{5}\right) =q^{6}+q^{4}\left( q-1\right) $
and $\ C_{1,2,3}=\left( 2q-1\right) ^{2}q^{3}$.

\ 

Then it is easy to see that $\left( 5.13\right) $\ is equivalent to $\left(
5.11\right) .$

\ 

Now assume that $a\neq 0.$ According to proposition IV-1 and $\left(
2.4a\right) $, we have:

\ 

\ $N\left( f,a\right) =$

$\quad \quad \sum\limits_{a_{1},a_{2}\in \mathbb{F}_{q}^{\ast
}}\sum\limits_{X_{6},X_{7},X_{8}}\prod\limits_{t=1}^{3}q^{\left| J_{\nu
_{t}}\right| -1}\left( 1+\kappa \left(
\sum\limits_{j=1}^{3}a_{t,j}^{^{\prime }}a_{j}-\sigma _{\left[ 1,2,3\right]
}^{t,1}X_{6}X_{7}X_{8}\right) \left( 1-\left( \frac{q-1}{q}\right) ^{\left|
J_{\nu _{t}}\right| -1}\right) \right) ,$

\ 

\ \ \ \ $B_{\left[ 1,2,3\right] }^{-1}=\left( a_{t,j}^{^{\prime }}\right)
_{t,j}=\left( 
\begin{array}{ccc}
1 & 0 & 0 \\ 
-1 & 1 & 0 \\ 
1 & -1 & 1%
\end{array}
\right) ,$ $B_{\left[ 1,2,3\right] }^{-1}C_{\left[ 1,2,3\right] }=\left( 
\begin{array}{c}
1 \\ 
0 \\ 
0%
\end{array}
\right) =\left( \sigma _{\left[ 1,2,3\right] }^{i,1}\right) _{i,1}.$

A straightforward computation gives us formula $\left( 5.12\right)
.\blacksquare $

\ 

\ 

VI \ WEIGHT HIERARCHY OF A CLASS OF CODES\medskip

Let us denote by $\mathcal{J}$ the partition $\{J_{1},J_{2},...,J_{m}\}$ of $%
\left\{ 1,...,n\right\} $. Following\ Cherdieu and Rolland in $[2]$, we can
introduce the map 
\begin{equation*}
\begin{array}{cccc}
c: & E(q,n,\mathcal{J}) & \rightarrow & \mathbb{F}_{q}^{q^{n}} \\ 
& f & \mapsto & \left( f\left( x\right) \right) _{x\in \mathbb{F}_{q}^{n}}%
\end{array}%
\end{equation*}
where $E(q,n,\mathcal{J})$ is the set of all multilinear polynomials with
separated variables of the form \ $f\left( X_{1},..,X_{n}\right)
=\sum_{j=1}^{m}a_{j}\prod_{\tau \in J_{j}}X_{\tau }$ where $%
a_{1},...,a_{m}\in \mathbb{F}_{q}$. The map $c$ is injective and its image $%
\func{Im}g$ is known as a multilinear code with separated variables. We
write $\func{Im}g=C\left( q,n,\mathcal{J}\right) $. Proposition III.1 gives: 
\begin{equation*}
N\left( f,a\right) =q^{n-1}+\left( q-1\right) q^{n-1}\left[
\prod\limits_{\left\{ j/a_{j}\neq 0\right\} }\left( 1-\left( \frac{q-1}{q}%
\right) ^{\left| J_{j}\right| -1}\right) \right] ,
\end{equation*}
and we easily deduce:\medskip

\textbf{PROPOSITION VI-1}

The\ lenght\ of\ the\ code\ $C\left( q,n,\mathcal{J}\right) $\thinspace
is\thinspace $q^{n}$,\thinspace and\thinspace its\thinspace dimension is $m$%
. The weight $\omega \left( c\left( f\right) \right) $ of the code-word $%
c\left( f\right) $, where $f\left( X_{1},..,X_{n}\right)
=\sum_{j=1}^{m}a_{j}\prod_{\tau \in J_{j}}X_{\tau }$, is 
\begin{equation*}
\omega \left( c\left( f\right) \right) =q^{n-1}\left( q-1\right) \left[
1-\prod\limits_{\left\{ j/a_{j}\neq 0\right\} }\left( 1-\left( \frac{q-1}{q}%
\right) ^{\left| J_{j}\right| -1}\right) \right] ,
\end{equation*}
and the minimum distance of $C\left( q,n,\mathcal{J}\right) $ is 
\begin{equation*}
\limfunc{dist}\left( C\left( q,n,\mathcal{J}\right) \right) =q^{n-s}\left(
q-1\right) ^{s}
\end{equation*}
where $s=\func{Max}_{1\leq j\leq m}\left( \left| J_{j}\right| \right) $%
.\medskip \medskip

If $C$ denotes a code with parameters $\left[ n,r,d\right] $, and if $D$
denotes one of its subcodes, we recall that $\chi \left( D\right) =\left\{
i\in \left\{ 1,...,n\right\} \,/\,\exists \left( x_{1},\ldots ,x_{n}\right)
\in D\ \text{with }x_{i}\neq 0\right\} $ is the support of $D$, and that the
number $\omega \left( D\right) =\left| \chi \left( D\right) \right| $ of
elements in $\chi \left( D\right) $ is the Wei weight of $D$ ($[6]$
p.~1412). The $h^{\text{th}}$\ minimum weight of $C$ is given by 
\begin{equation*}
d_{h}=d_{h}\left( C\right) =\func{Min}\left\{ \omega \left( D\right) \,/\,D%
\text{ subcode of }C\text{ with }\dim (D)=h\right\} .
\end{equation*}
The minimum distance of $C$ is $d_{1}$, and the weight hierarchy of the code 
$C$ is the sequence $\left\{ d_{1},\ldots ,d_{r}\right\} .$\medskip

If $D=\limfunc{Vect}\left( f_{1},...,f_{h}\right) $ denotes the subcode of $%
C\left( q,n,\mathcal{J}\right) $ generated by the polynomials $f_{1}$, ..., $%
f_{h}$ in $E(q,n,\mathcal{J})$, and if we write $\mathbb{F}_{q}^{n}=\left\{
P_{1},...,P_{q^{n}}\right\} $, the Wei weight of $D$ is 
\begin{eqnarray*}
\omega \left( D\right) &=&\left| \left\{ i\in \left\{ 1,...,q^{n}\right\}
\;/\;\exists g\in D\text{\quad }g\left( P_{i}\right) \neq 0\right\} \right|
\\
&=&\left| \left\{ i\in \left\{ 1,...,q^{n}\right\} \;/\;\exists \lambda
_{1},...,\lambda _{h}\in \mathbb{F}_{q}\text{\quad }\sum\limits_{j=1}^{h}%
\lambda _{j}f_{j}\left( P_{i}\right) \neq 0\right\} \right| \\
&=&\left| \left\{ i\in \left\{ 1,...,q^{n}\right\} \;/\;\exists j\in \left\{
1,...,h\right\} \ \text{\quad }f_{j}\left( P_{i}\right) \neq 0\right\}
\right| ,
\end{eqnarray*}
thus 
\begin{equation*}
\omega \left( D\right) =q^{n}-N\left( f_{1},...,f_{h},0,...,0\right) .
\end{equation*}
Each polynomial $f_{i}$ ($1\leq i\leq h$) can be written in the form $\
f_{i}\left( X_{1},..,X_{n}\right) =\sum_{j=1}^{m}a_{ij}\prod_{\tau \in
J_{j}}X_{\tau }$ and the rank of the matrix $A=\left( a_{ij}\right) $ is $k$
as soon as we assume that $D$ is of dimension $k$. We can therefore apply
Proposition IV.1 to compute $\omega \left( D\right) $ and $d_{h}$.\bigskip
\bigskip\ 

VII \ \ CONCLUSION

\ 

As we have seen, Proposition IV.1 gives us a semi-explicit formula to
compute the number $N\left( f_{1},...,f_{k},0,...,0\right) $ when the rank
of the matrix $A$ is $k$, and it would be interesting to build an explicit
program to do it, and compute the weight hierarchy of the code $C\left( q,n,%
\mathcal{J}\right) $. This could be a numerical continuation of this
work.\bigskip

Several problems can arise from the present paper, and we gives here two
ways for further investigations. First, it would be interesting to define a
partition $\left\{ \triangle _{1},\triangle _{2},...,\triangle _{\eta
}\right\} $ of $\mathbb{F}_{q}^{n}$ and write the number $N\left( f,a\right) 
$ in this way: 
\begin{equation*}
N\left( f,a\right) =\sum_{j=1}^{\eta }N_{\triangle _{j}}\left( f,a\right)
=\sum_{j=1}^{\eta }\sum_{a_{1}...a_{k}=a}N_{\triangle _{j}}\left(
f_{1},...,f_{k},a_{1},...,a_{k}\right) ,
\end{equation*}
where $N_{\triangle _{j}}\left( f,a\right) $ denotes the number of
solutions\ of $\left( 1.3\right) $ in $\triangle _{j}$, and $N_{\triangle
_{j}}\left( f_{1},...,f_{k},a_{1},...,a_{k}\right) $ stands for the number
of solutions\ of the system $f_{i}\left( X_{1},...,X_{n}\right) =a_{i}$ ($%
1\leq i\leq k$) in $\triangle _{j}$. Specific partitions could lead us to
less calculus. For instance,\medskip

\quad $\triangle _{1}=\left\{ \left( u_{1},...,u_{n}\right) \in \mathbb{F}%
_{q}^{n}\,/\text{\thinspace }\exists \left( i_{1},i_{2},..,i_{m}\right) \in
J_{1}\times J_{2}\times ...\times J_{m}\quad \left(
u_{i_{1}},...,u_{i_{m}}\right) =0\right\} ,\smallskip $

\quad $\triangle _{2}=\left\{ \left( u_{1},...,u_{n}\right) \in \left( 
\mathbb{F}_{q}^{\ast }\right) ^{n}\,/\text{\thinspace }\forall \tau \in
\left\{ 1,..,m\right\} \quad \forall \left( i,j\right) \in J_{\tau }\times
J_{\tau }\quad u_{i}=u_{j}\ \text{for all }\right\} ,$\medskip

lead to a trivial numbers and diagonal equations.\medskip

The second problem consists to study the more general case where:\medskip

\quad $f=f_{1}f_{2}...f_{k},$

\quad $f_{i}\left( X_{1},...,X_{n}\right) =\sum\limits_{j=1}^{s}$\ $%
a_{ij}\prod\limits_{\tau \in J_{j}^{\left( i\right) }}X_{\tau }$, $i=1,..,k,$

\quad $\left\{ J_{1}^{\left( 1\right) },J_{2}^{\left( 1\right)
},...,J_{s}^{\left( 1\right) }\right\} $, $\left\{ J_{1}^{\left( 2\right)
},J_{2}^{\left( 2\right) },...,J_{s}^{\left( 2\right) }\right\} $, ..., $%
\left\{ J_{1}^{\left( k\right) },J_{2}^{\left( k\right) },...,J_{s}^{\left(
k\right) }\right\} $ are partitions of $\left\{ 1,...,n\right\} $,

and to obtain semi-explicit formulas for the number of solutions of $f\left(
X_{1},...,X_{n}\right) =a$ in this case.

\ 

\ 

REFERENCES

\ \ 

$\left[ 1\right] $ \ L. Carlitz, The number of solutions of some special
equations in a finite field, Pacific J. Math.~4, 1954, pp.~207-217.

\ 

$\left[ 2\right] $ \ J.-P. Cherdieu and\ R. Rolland, On hypersurfaces
defined by a separated variables polynomial over a finite field, Arithmetic,
geometry and coding theory (Luminy, 1993), De Gruyter, Berlin, 1996,
pp.~35-43.

\ 

$\left[ 3\right] $ \ J.-R. Joly, Equations et vari\'{e}t\'{e}s alg\'{e}%
briques sur un corps fini, Enseignement Math. 19, 1973, pp.~1-117.

\ 

$\left[ 4\right] $ \ A. K. Lenstra, Factoring multivariate polynomials over
finite fields, Journal of computer and system sciences 30, 1985, pp.~235-248.

\ 

$\left[ 5\right] $ \ R.G. Van Meter, The number of solutions of certain
systems of equations in a finite field,\ Duke\ Math. Journal 38, 1971,
pp.~365-377.

\ 

$\left[ 6\right] $ \ V. Wei, Generalized hamming weights for linear codes,
IEEE Transactions on information theory 37, 1991, pp.~1412-1418.

\end{document}